\documentclass[reqno,11pt]{amsart}
\usepackage{amsmath, amssymb, amsthm,amsfonts} 
\usepackage[english]{babel}
\usepackage{bbm}
\usepackage{graphicx}
\usepackage{xcolor}
\usepackage{url}
\usepackage{epstopdf}
\usepackage[a4paper,bindingoffset=0.5cm,left=2cm,right=2cm,top=2.5cm,bottom=2cm,footskip=.8cm]{geometry}

\usepackage{tikz}
\usetikzlibrary{arrows,automata,positioning,calc,shapes,decorations.pathreplacing,decorations.markings,shapes.misc,petri,topaths}
\usepackage{tkz-berge}
  \usepackage{pgfplots}
  \pgfplotsset{compat=newest}
  \usetikzlibrary{plotmarks}
\newlength\figureheight
  \newlength\figurewidth
\pgfplotsset{%
    tick label style={font=\scriptsize},
    label style={font=\footnotesize},
    legend style={font=\footnotesize},
         every axis plot/.append style={very thick}
}

\usepackage{rotating}
\usepackage{amsbsy,enumerate}
\usepackage{graphicx}
\usepackage{ccaption}
\usepackage{comment}
\usepackage{mathrsfs}

\newcommand{\vb}{\vspace{3.2mm}}
\renewcommand{\hat}{\widehat}

\makeatletter
\renewcommand{\fnum@figure}[1]{\textbf{\figurename~\thefigure}. }
\renewcommand{\fnum@table}[1]{\textbf{\tablename~\thetable}. }
\makeatother

\usepackage{amsmath, amssymb, amsfonts, amsthm, bm, mathtools, dsfont}
\usepackage{tikz}

\let\originalleft\left
\let\originalright\right
\renewcommand{\left}{\mathopen{}\mathclose\bgroup\originalleft}
\renewcommand{\right}{\aftergroup\egroup\originalright}

\newcommand{\abs}[1]{\left\vert #1 \right\vert}
\newcommand{\absj}[1]{| #1 |}
\newcommand{\diag}[1]{\mathrm{diag} \left( #1 \right)}
\newcommand{\ind}[1]{\mathds{1}_{\left\lbrace #1 \right\rbrace}}

\newcommand{\norm}[1]{\left\| #1 \right\|}
\newcommand{\transpose}{\mathsf{T}}

\newcommand{\dd}{\mathrm{d}}

\newcommand{\genQ}{\mathcal{Q}}

\newcommand{\prequadvar}[1]{\left\langle #1 \right\rangle}
\newcommand{\stdot}{\mathbin{\bm{\cdot}}}

\newcommand{\eps}{\epsilon}

\theoremstyle{plain}

\newtheorem{theorem}{Theorem}[section]
\newtheorem{assumption}[theorem]{Assumption}
\newtheorem{corollary}[theorem]{Corollary}
\newtheorem{lemma}[theorem]{Lemma}

\theoremstyle{definition}

\newtheorem{example}[theorem]{Example}
\newtheorem{remark}[theorem]{Remark}

\allowdisplaybreaks

\title[Networks of {M}arkov-modulated  infinite-server queues]{Diffusion limits for networks of {M}arkov-modulated  infinite-server queues}
\author{H. M. Jansen$^{1,2}$, M.  Mandjes$^{1}$, K. De Turck$^{3}$, S. Wittevrongel$^{2}$}

\begin{document}
\maketitle

\begin{abstract}
\noindent This paper studies the diffusion limit for a network of infinite-server queues operating under Markov modulation (meaning that the system's parameters depend on an autonomously evolving background process). In previous papers on (primarily single-node) queues with Markov modulation, two variants were distinguished: one in which the server speed is modulated, and one in which the service requirement is modulated (i.e., depends on the state of the background process upon arrival). The setup of the present paper, however, is more general, as we allow both the server speed and the service requirement to depend on the background process. For this model we derive a Functional Central Limit Theorem:
we show that, after accelerating the arrival processes and the background process, a centered and normalized version of the network population vector  converges to a multivariate Ornstein-Uhlenbeck process. The proof of this result relies on expressing the queueing process in terms of Poisson processes with a random time change, an application of the Martingale Central Limit Theorem, and continuous-mapping arguments. 

\vb


\vb

\noindent $^{1}$ Korteweg-de Vries Institute for Mathematics,
University of Amsterdam, Science Park 904, 1098 XH Amsterdam, the Netherlands.

\noindent $^{2}$ TELIN, Ghent University, Sint-Pietersnieuwstraat 41,
B-9000 Ghent, Belgium.

\noindent $^{3}$ CentraleSup\'elec,
D\'epartement de T\'el\'ecommunications,
Laboratory of Signals and Systems (L2S),
UMR8506
Plateau de Moulon,
3 rue Joliot-Curie,
91192 Gif sur Yvette,
France.

\noindent {\it E-mail}. {\tt\{h.m.jansen|m.r.h.mandjes\}@uva.nl}, {\tt koen.deturck@supelec.fr}, {\tt sw@telin.ugent.be}

\vb

\noindent
M.\ Mandjes' research is partly funded by the NWO Gravitation Programme N{\sc etworks}, Grant Number 024.002.003, and an NWO Top Grant, Grant Number 613.001.352.

\end{abstract}

%

\section{Introduction}
\label{sec:intro}
Networks of infinite-server queues with time-dependent parameters have been extensively studied. In a series of articles, several nice properties have been shown for the situation that the network's queues are fed by Poisson arrival streams with time-dependent deterministic arrival rates. In this model, due to the fact that customer streams do not interact, departure streams are again time-dependent Poisson (under specific assumptions on the queues' initial conditions). In addition, at any point in time the  joint queue length has independent Poisson marginals. Without pursuing an exhaustive overview, we refer to the contributions by Harrison and Lemoine~\cite{HarLem} and Eick {\it et al.}~\cite{Eick}. 

Modelling customer streams as Poisson processes with time-dependent deterministic arrival rates has, however, the following crucial drawback. Due to standard properties of such processes, for any time interval the mean number of arrivals coincides with the corresponding variance. In various practical contexts, however, there is {\it overdispersion}, i.e., the variance substantially exceeds the mean; see e.g.\ \cite{Zeevi,KimWhi} and references therein. 
To  incorporate overdispersion a common procedure is to apply {\it modulation}, meaning that the arrival rate depends on an independently evolving {stochastic} process (usually referred to as the {background process}). Due to the additional uncertainty that is introduced by the background process, the variance of the number of arrivals now exceeds the mean. An example of modulation is to let the arrival process be a Cox process, i.e., a Poisson process whose rate is a (nonnegative) stochastic process. 

Above we introduced the concept of modulation for arrival processes, but it can clearly be applied to the queues' service processes as well. 
It is important to notice that modulation of service processes comes in multiple flavors. In the first place, we can modulate the server speed. A slightly more complicated form of modulation lets the service requirement distribution depend on the background process: the service requirement is sampled upon arrival, with a parameter determined by the state of the background process at the moment the customer arrives. Thus, based on which parameters of the model are modulated, the following four types of modulation may be distinguished. In line with the terminology used in \cite{BDM2016a},
\begin{itemize}
\item[$\circ$]
\makebox[1.8cm]{Model~0:\hfill} only the arrival rate is modulated,
\item[$\circ$]
\makebox[1.8cm]{Model~I:\hfill} only the arrival rate and the server speed are modulated,
\item[$\circ$]
\makebox[1.8cm]{Model~II:\hfill} only the arrival rate and the service requirement distribution are modulated,
\item[$\circ$]
\makebox[1.8cm]{Model~III:\hfill} all parameters are modulated.
\end{itemize}
Obviously, Model~0 is a special case of both Model~I and Model~II, and Model~I and Model~II are special cases of Model~III.

Introducing modulation complicates the analysis of networks of infinite-server queues considerably: many of the nice properties that were found for time-dependent deterministic arrival rates do not carry over to the setting in which the arrival rate is stochastic, for instance. More specifically, the joint queue length ceases to have independent marginals, and in addition these marginals are not Poisson. An important 
strand of work is focused on {\it Markov modulation}: the value of the model parameters is determined by an external finite-state continuous-time Markov chain. Early work on (primarily single-node) Markov-modulated infinite-server systems, such as~\cite{OP1986}, develops procedures to identify the (transient and stationary) moments of the number of clients; see also \cite{BX2004,Dauria2008, FA2009} for more recent contributions.  A significant drawback of most of these results concerns their implicit nature: systems of partial differential equations are derived for the moment generating function pertaining to the joint distribution of the number of clients in the individual queues, and recursions for the corresponding moments are found. Only in very special cases the underlying probability distribution can be found in more explicit terms; see e.g.\ \cite{Falin2008}.

The lack of explicit expressions motivated the systematic study of scaling regimes in which the queue length distribution {\it can} be given in closed form. Some recent papers have established scaling limits in the context of single-queue Model~0 systems: \cite{HLM2017} considers the scenario in which the arrival rate is periodically resampled, whereas \cite{KBM2017} analyzes the case in which the arrival rate is given by a shot-noise process. For the case of single queues under Markov modulation, it is shown in \cite{ABMTD2016} (for Model~0) and \cite{BDM2016a} (for Model~I and Model~II) that, after appropriately accelerating the arrival processes and the background process, a centered and normalized version of the number of customers in the system converges to an Ornstein-Uhlenbeck process. The proof of this Functional Central Limit Theorem (FCLT) relies on casting the queueing process in terms of Poisson processes with a random time change and an application of the Martingale Central Limit Theorem.

The main goal of the present paper is to extend the findings of \cite{ABMTD2016,BDM2016a} in two directions. First of all, we study a network of infinite-server queues under Markov modulation (rather than a single queue). Next to that, we do not restrict ourselves to Model~0, Model~I, or Model~II, but focus on the Model~III variant which covers the other variants as special cases. Our contribution is that we derive, for this highly general setting, an FCLT. 

The different setting that we study also motivates the use of a different strategy to establish an FCLT. As will become clear in the paper, the argumentation for this multivariate setting has some elements in common with the one used in the single-queue case, but in general requires considerably more delicate arguments. A key insight is that we may exploit continuous-mapping arguments once we have established asymptotic equivalence of our network with some other, less complicated stochastic process. Following this strategy, we succeed in establishing an FCLT: the centered and normalized version of the joint queue length process converges to a multivariate Ornstein-Uhlenbeck process. 

The paper is organized as follows. In Section \ref{sec:model}, we introduce our model and point out how to express the network population vector in terms of unit-rate Poisson processes with a random time change. We also show that, perhaps somewhat surprisingly, any network of Model III type can be expressed in terms of a, seemingly more restrictive, Model I network, entailing that it suffices to just consider the Model I setting from then on. 
Then Section~\ref{sec:fluidlimit} presents the fluid limit under a scaling of both the arrival rates and the background process. Centering the queue lengths with this fluid limit, the results in Section \ref{sec:diffusionlimit} show that, after a suitable normalization, the joint queue length process weakly converges to a multivariate Ornstein-Uhlenbeck process. 
Proofs and auxiliary results are given in Section~\ref{proofs}.

\section{Model and preliminaries}
\label{sec:model}
In this section we define a Markov-modulated network of infinite-server queues with exponentially distributed service requirements. We introduce notation for the Model I variant and 
show how the queue length vector can be represented in terms of Poisson processes with a random time change. We then argue that, in the setting that we consider, a network containing Model II or Model III queues can be represented as a network containing only Model I queues (justifying that we restrict ourselves to Model I in the rest of the paper). 

\subsection{Definition of Model I}
We now provide the mathematical description of the Model I version of a Markov-modulated network of infinite-server queues. We assume that the background process $J$ is an autonomously evolving continuous-time Markov chain with finite state space $\left\lbrace 1 , \dotsc , d \right\rbrace$ and irreducible transition rate matrix $\genQ$. The invariant distribution corresponding to $\genQ$ is denoted by the $d \times 1$ column vector $\pi$ (so $\pi$ has nonnegative entries that sum to $1$ and $\pi^{\transpose} \genQ = 0$). 

Jobs that leave a particular queue may be routed to another queue or leave the network altogether. 
Throughout the paper we assume there are $L \in \{2,3,\ldots\}$ queues in the network.
We impose the assumption of at least two queues, because this allows us to not explicitly record the jobs leaving the system. Indeed, a job leaving the system is equivalent to sending it to some designated infinite-server queue that has zero server speed (i.e., a queue from which jobs cannot leave). The assumption $L\in\{2,3,\ldots\}$ is therefore not a restriction and is merely intended to streamline notation.

Given a network of $L$ queues, we define the set of queue indices $\mathcal{I} = \left\lbrace 1 , \dotsc , L \right\rbrace$. For any $k \in \mathcal{I}$, we define $\mathcal{I}_{k} = \mathcal{I} \setminus \left\lbrace k \right\rbrace$. This is, of course, the set of all queues to which a job from queue $k$ may be routed after service completion. This set is always nonempty, because $L \in \{2,3,\ldots\}$ by assumption.

When the background process is in $i\in\{1,\ldots,d\}$, jobs arrive at queue $k$ as a Poisson process with rate $\lambda_k(i)$. In addition, again for the background process being in state $i$, we denote  by $\mu_{kl} \left( i \right)$ the rate at which any individual job present in queue $k$ completes service and jumps to queue $l$. Observe that during the time a job spends in a given queue, the rate of jumping to another queue may change in time (when the state of the background process changes, that is).

We denote the number of jobs in queue $k \in \mathcal{I}$ at time $t\geq 0$ by $Q_{k} \left( t \right)$, and the $L$-dimensional queue content process by
\begin{align}
Q \left( t \right) = \left( Q_{1} \left( t \right) , \dotsc , Q_{L} \left( t \right) \right).
\label{eq:queuecontentprocess}
\end{align}
Then we represent $Q$ as the solution of the system of equations
\begin{align}
\begin{split}
Q_{k} \left( t \right)
&= 
Q_{k} \left( 0 \right) + A_{k} \left( \int_{0}^{t} \lambda_{k} \left( J \left( s \right) \right) \, \dd s \right) + \sum_{l \in \mathcal{I}_{k}} S_{lk} \left( \int_{0}^{t} \mu_{lk} \left( J \left( s \right) \right) Q_{l} \left( s \right) \, \dd s \right) \\
&\phantom{{} = {}} {} - \sum_{l \in \mathcal{I}_{k}} S_{kl} \left( \int_{0}^{t} \mu_{kl} \left( J \left( s \right) \right) Q_{k} \left( s \right) \, \dd s \right).
\end{split} \label{eq:Poissonprocessequations}
\end{align}
Here, $Q \left( 0 \right)$ is an independent random vector taking values in $\mathbb{Z}_{\geq 0}^{L}$, with its $k$-th entry describing the initial number of jobs in queue $k$. The processes $A_{k}$ and $S_{kl}$ with $k \in \mathcal{I}$ and $l \in \mathcal{I}_{k}$ form a collection of independent standard (unit-rate, that is) Poisson processes. We use $A_{k}$ to model the exogenous arrivals to queue $k$, whereas we use $S_{kl}$ to model service completions in queue $k$ that are routed to queue $l$. This  representation has been extensively used in e.g.\ \cite{PTW2007}, and is (in our experience at least) an intuitive representation that is relatively easy to work with when studying weak convergence.

\subsection{Models II and III are covered by Model I} In the network setting that we consider, it turns out, perhaps somewhat surprisingly, that any Model II or Model III network of infinite-server queues can be described through a Model I network of infinite-server queues.  At first sight, our construction of a network of Model I queues may not seem very flexible, because the modulation of service requirements (in Model II and Model III queues) is not explicitly incorporated. However, the combination of Model I queues with the network setting gives us much more flexibility than we might have thought.

The upcoming examples are intended to demonstrate this flexibility by explicitly showing how we may incorporate modulated service requirements. In particular, the second example demonstrates how we may construct a Model III queue as a network of Model I queues, thus showing that a network of Model I queues is indeed the most general Markov-modulated network of infinite-server queues.

\begin{remark}\label{REM1}
Observe that the Model I setting incorporates probabilistic routing: when the background process is in state $i$, the rate of leaving queue $k$ is $\mu_k(i):=\sum_{l\not = k}\mu_{lk}(i)$, and the probability that after being served the job jumps to queue $l$ is $p_{kl}(i):=\mu_{kl}(i)/\mu_k(i).$ As such, the routing probabilities depend on the state of the background process at service completion. $\hfill\Diamond$
\end{remark}

%
%

\begin{example}
In Remark \ref{REM1} we already argued that Model I covers settings in which the routing probabilities depend on the current state of the background process. An obvious extension would be to make the routing probabilities depend on the state of the background process at the moment of a job's arrival. In this example we argue that such a mechanism still falls in the context of Model I networks of Markov-modulated infinite-server queues.

We consider a network in which jobs arrive to queue $\textsc{a}$ at rate $\lambda_{\textsc{a}} \left( J \left( t \right) \right)$ and receive service there with speed $\mu_{\textsc{a}} \left( J \left( t \right) \right)$. When a job completes service at this queue, it is routed with probability $p \left( i , j \right)$ to queue $\textsc{b}$, where $i$ is the state of the background process $J$ upon arrival of the job and $j$ is the state of $J$ upon service completion. With probability $1 - p \left( i , j \right)$ the job is routed to queue $\textsc{c}$. (What happens at queue $\textsc{b}$ and $\textsc{c}$ is not very relevant here; neither do we need a queue $\textsc{d}$ for jobs leaving the system.)
The difficulty here is the dependence of the routing probabilities on both the present and the past. Indeed, the routing probabilities for a job are determined when it completes service, but also depend on the state of $J$ upon its arrival.

We can deal with this problem by registering the state of the background process $J$ when a job arrives to queue $\textsc{a}$. This naturally leads to the jobs in queue $\textsc{a}$ being divided into $d$ different classes (where we recall that $d$ is the size of the state space of $J$): a job belongs to class $i$ when $J$ is in state $i$ upon arrival of the job. The idea is now to split the content of queue $\textsc{a}$ into  $d$ individual queues, i.e., one for each class. Then the routing probabilities for jobs in those queues only depend on the queue index and the state of $J$ upon service completion, which is in line with the probabilistic routing mechanism that we discussed in Remark \ref{REM1}.

In more concrete terms, the system considered as a network of Model I queues will look like this. We denote the queues handling the different classes by $1 , \dotsc , d$. We define the arrival rates to these queues by $\lambda_{k} \left( i \right) = \ind{i=k} \lambda_{\textsc{a}} \left( i \right)$ for $k,i \in \left\lbrace 1 , \dotsc , d \right\rbrace$. Then we have arrivals of only one class at a time with the right arrival rate for each class. 

Jobs from those $d$ queues will be routed to queue $d+1$ and queue $d+2$ (which were queue $\textsc{b}$ respectively queue $\textsc{c}$ in our original thought experiment). The routing mechanism is then incorporated by defining the server speeds by picking $\mu_{k \left( d+1 \right)} \left( i \right) = p \left( k , i \right) \mu_{\textsc{a}} \left( i \right)$ and $\mu_{k \left( d+2 \right)} \left( i \right) = \left( 1 - p \left( k , i \right) \right) \mu_{\textsc{a}} \left( i \right)$ for $k \in \left\lbrace 1 , \dotsc , d \right\rbrace$. $\hfill\Diamond$
\end{example}

\begin{example}
In the previous example, we have seen how we can handle different classes (corresponding to the state of the background process that the job sees upon arrival) by assigning them to their own Model I queue. In this example, we use this idea to show how a Model III Markov-modulated infinite-server queue may be represented as a Model I Markov-modulated network of  infinite-server queues.

Consider a single Model III Markov-modulated infinite-server queue where $J$ is the background process. In this case, jobs arrive according to a Poisson process with rate $\lambda^{*} \left( i \right)$ while $J$ is in state $i$. If a job arrives to the system while $J$ is in state $i$, then it will be of type $i$ and has an independent service requirement having an exponential distribution with parameter $\kappa^{*} \left( i \right)$. The server speed is $\mu^{*} \left( i \right)$ while  the background process is in state $i$.

We will construct a network of $L = d+1$ queues: $d$ queues for each job type (i.e., queue $k$ contains the jobs of type $k$ that are still in service) and one queue to collect all jobs that have finished service in their respective queue (i.e., queue $d+1$ with zero exogenous arrival rate and zero server speed). In what follows, $i$ ranges over $\left\lbrace 1 , \dotsc , d \right\rbrace$, and $k$ and $l$ range over $\left\lbrace 1 , \dotsc , d , d+1 \right\rbrace$.

To construct the network, define $\lambda_{k} \left( i \right) = \ind{i = k} \lambda^{*} \left( i \right)$ for $k \in \left\lbrace 1 , \dotsc , d \right\rbrace$ and $\lambda_{d+1} \left( i \right) = 0$. This means that jobs arrive to exactly one queue at a time: the arrival rate (of jobs of type $i$) to queue $i$ is $\lambda^{*} \left( i \right)$ while $J$ is in state $i$ and zero while $J$ is not in state $i$. Moreover, there are no exogenous arrivals to queue $d+1$, which collects jobs that have finished service.
Additionally, we take $\mu_{kl} \left( i \right) = \ind{l = d+1} \kappa^{*} \left( k \right) \mu^{*} \left( i \right)$ for $k \in \left\lbrace 1 , \dotsc , d \right\rbrace$. This ensures that jobs can only be routed to queue $d+1$ after service completion (the $\ind{l = d+1}$ part), that jobs of type $k$ arriving at queue $k$ have an exponentially distributed service requirement with parameter $\kappa^{*} \left( k \right)$, and that those jobs in queue $k$ experience server speed $\mu^{*} \left( i \right)$ while $J$ is in state $i$. 

We have thus defined a Model I network of $d+1$ infinite-server queues, $Q_1(t)$ up to $Q_{d+1}(t)$, that is equivalent to our single Model III Markov-modulated infinite-server queue. 
Now the queue content process of this single Model III queue at time is  given by $\sum_{k=1}^{d} Q_{k}(t)$. As a bonus, the departure process of the Model III queue is captured by $Q_{d+1}(t)$. Thus, the Model I network constructed above completely describes the behavior of the queue content of the Model III queue.

Along the same lines, it can be argued that a network of $L$ Model III Markov-modulated infinite-server queues can be described as a network of $Ld$ Model I Markov-modulated infinite-server queues. We conclude that we can restrict ourselves to Model I.
$\hfill\Diamond$
\end{example}

\begin{remark}
It is directly seen that phase-type service requirements are also easy to deal with, by introducing infinite-server queues corresponding to each of the phases. $\hfill\Diamond$
\end{remark}

\section{Convergence results: a fluid limit}
\label{sec:fluidlimit}
As mentioned earlier, we are interested in proving a diffusion limit for our Markov-modulated network of infinite-server queues. A first step towards such a result is the derivation of a fluid limit for the network. This serves three purposes. First, the fluid limit is interesting in its own right, as it describes the typical behavior of the network under scaling. It is, essentially, a law of large numbers for the evolution of the network population. It holds under more general assumptions than the diffusion limit. Second, in the proof for the diffusion limit we will apply the fluid limit to identify the limiting process; the fluid limit limit is also used to appropriately center the network population process. Third, we will use the fluid limit to point out differences between the modulated network and its nonmodulated counterpart.

To establish the fluid limit result, we scale the parameters of the modulated network. We do so by scaling the arrival rates to the queues, the server speeds, and the transition rates of the background process. In addition, the network's initial population is scaled appropriately. More specifically, we impose the following assumption.

\begin{assumption}\label{AS1}
In the $n$-th system,
\begin{itemize}
\item[$\circ$]
with $Q_{k}^{n} \left( 0 \right)$ denoting the initial number of jobs in queue $k$, we assume  that $\tfrac{1}{n} Q^{n} \left( 0 \right)$ converges in probability to the constant $\rho \left( 0 \right)$;
\item[$\circ$] with $\lambda_{k}^{n}$ the arrival rate to queue $k$, we assume that
\begin{align}
\lim_{n \to \infty} \frac{1}{n} \lambda_{k}^{n} = \lambda_{k};
\label{eq:lambdakfluid}
\end{align}
\item[$\circ$] with $\mu_{kl}^{n}$
 the service speed for jobs to be routed from queue $k$ to queue $l$, we assume that\begin{align}
\lim_{n \to \infty} \mu_{kl}^{n} = \mu_{kl};
\label{eq:muklfluid}
\end{align}\item[$\circ$] 
the background process, denoted by $J^{n}$, corresponds to a continuous-time Markov chain having irreducible generator matrix $n^{\alpha} \genQ$ (for some $\alpha>0$) and invariant distribution $\pi$. \end{itemize}
\end{assumption}
Regarding the $\mu_{kl}^n$ and $\lambda^n_k$, the main examples to keep in mind are $\lambda_{k}^{n} = n \lambda_{k}$ and $\mu_{kl}^{n} = \mu_{kl}$.

Recalling that $J^{n} \left( t \right)$ matches with $J \left( n^{\alpha} t \right)$, we observe that the scaling proposed here amounts to simultaneously speeding up the arrivals and the time scale of the background process, while the server speeds approach some limiting value. Note that the speedup of the arrival rate is essentially linear in $n$, whereas the speedup of the time scale of the background process is sublinear ($\alpha < 1$), linear ($\alpha = 1$), or superlinear ($\alpha > 1$) in $n$.

Under this scaling, the queue content process is denoted by the $\mathbb{R}^{L}$-valued stochastic process $Q^{n}$, which is given by a system of equations involving unit-rate Poisson processes with a random time change:
for $k=1,\ldots,L$ we have
\begin{align}
\begin{split}
Q_{k}^{n} \left( t \right)
&= 
Q_{k}^{n} \left( 0 \right) + A_{k} \left( \int_{0}^{t} \lambda_{k}^{n} \left( J^{n} \left( s \right) \right) \, \dd s \right) + \sum_{l \in \mathcal{I}_{k}} S_{lk} \left( \int_{0}^{t} \mu_{lk}^{n} \left( J^{n} \left( s \right) \right) Q_{l}^{n} \left( s \right) \, \dd s \right) \\
&\phantom{{} = {}} {} - \sum_{l \in \mathcal{I}_{k}} S_{kl} \left( \int_{0}^{t} \mu_{kl}^{n} \left( J^{n} \left( s \right) \right) Q_{k}^{n} \left( s \right) \, \dd s \right).
\end{split} \label{eq:scaledsysdef}
\end{align}
(See also Eq.\ \ref{eq:Poissonprocessequations}.)
As mentioned, we are interested in finding a fluid limit for $Q^{n}$. This will involve the time-averaged parameters $\lambda_{k}^{\pi} = \sum_{i=1}^{d} \pi \left( i \right) \lambda_{k} \left( i \right)$ and $\mu_{kl}^{\pi} = \sum_{i=1}^{d} \pi \left( i \right) \mu_{kl} \left( i \right)$.
We present the fluid limit in the next lemma.
\begin{lemma}
\label{lem:fluidlimit}
Under Assumption \ref{AS1}, the process $\tfrac{1}{n} Q^{n}$ converges uoc in probability to the solution of the system of integral equations
\begin{align}
\rho_{k} \left( t \right)
&= \rho_{k} \left( 0 \right) + \int_{0}^{t} \lambda_{k}^{\pi} \, \dd s + \sum_{l \in \mathcal{I}_{k}} \int_{0}^{t} \mu_{lk}^{\pi} \rho_{l} \left( s \right) \, \dd s - \sum_{l \in \mathcal{I}_{k}} \int_{0}^{t} \mu_{kl}^{\pi} \rho_{k} \left( s \right) \, \dd s, \label{eq:rhodef}
\end{align}
for $k=1,\ldots,L$.
\end{lemma}
This fluid limit is exactly the fluid limit that we would obtain from a {\it nonmodulated} network of infinite-server queues, and is in particular not affected by the value of $\alpha$. (To see this, just replace the modulated parameters of the queues by their time-averaged versions $\lambda_{k}^{\pi}$ and $\mu_{kl}^{\pi}$.) However, modulation induces markedly different fluctuations {\it around} this fluid limit, i.e., in the diffusion scaling. We  explore this in the next section.

\section{Convergence results: a diffusion limit}
\label{sec:diffusionlimit}
Our main goal in this section is to establish a diffusion limit for the queue content process of the network under an appropriate scaling (which differs from the scaling imposed in the fluid limit result). This diffusion limit depends explicitly on scaling properties of the background process, thus giving insight into how the background process influences the behavior of the network on process level.

We will use several different tools to prove the diffusion limit. Typically, diffusion limits are proven relying on martingale methods and the Continuous Mapping Theorem (CMT); see, for instance, \cite{PTW2007,Whitt1980,whitt2002}. Unfortunately, the classical continuous maps do not apply to our network, due to the presence of the background process. We may get around this problem by constructing a more refined continuous map, but we opt for another solution. A crucial step in this solution consists of showing the asymptotic equivalence of the scaled network queue content process and a closely related stochastic process, which turns out to have the desired weak convergence properties. This approach combines queueing intuition with well-known continuous maps and convergence properties of stochastic integrals to achieve these results.

The continuous map that we will use is (a multidimensional version of) the well-known integral map presented in \cite[Th.~4.1]{PTW2007}. It is intimately related to standard infinite-server queues (cf.\ \cite{PTW2007} and Lemma~\ref{lem:fluidlimit}). For completeness, we present the result here. Its proof is a straightforward generalization of the proof of \cite[Th.~4.1]{PTW2007}.
\begin{lemma}
\label{lem:ctuintegralmap}
Let $h \colon \mathbb{R}^{L} \to \mathbb{R}^{L}$ be Lipschitz continuous, i.e.,
\begin{align*}
\norm{h \left( a_{1} \right) - h \left( a_{2} \right)} \leq c \norm{ a_{1} - a_{2} }
\end{align*}
for all $a_{1}, a_{2} \in \mathbb{R}^{L}$ for some $c > 0$. Also assume that $h \left( 0 \right) = 0$. If $b \in \mathbb{R}^{L}$ and $x \in D \left( \left[ 0,\infty \right) ; \mathbb{R}^{L} \right)$, then there exists a unique element $y \in D \left( \left[ 0,\infty \right) ; \mathbb{R}^{L} \right)$ satisfying the integral equation
\begin{align}
y \left( t \right) &= b + x \left( t \right) + \int_{0}^{t} h \left( y \left( s \right) \right) \, \dd s. \label{eq:ctuintegralmap}
\end{align}
Thus, the integral equation \eqref{eq:ctuintegralmap} defines a function $H \colon \mathbb{R}^{L} \times D \left( \left[ 0,\infty \right) ; \mathbb{R}^{L} \right) \to D \left( \left[ 0,\infty \right) ; \mathbb{R}^{L} \right)$. The function $H$ is continuous whenever the space $D \left( \left[ 0,\infty \right) ; \mathbb{R}^{L} \right)$ is equipped (in both the domain and the range) either with the uniform topology or with the weak J$_{1}$ topology.
\end{lemma}

The next corollary describes an example of a function $H$ that is defined through the procedure in Lemma~\ref{lem:ctuintegralmap}. The function $H$ defined in this corollary is, in fact, the map that we will use to establish the diffusion limit. Moreover, $H$ is clearly related to the fluid limit $\rho$ (to see this, compare Eq.\ \eqref{eq:rhodef} and Eq.\ \eqref{eq:ctuintegralmapnetwork}).

\begin{corollary}
\label{ex:ctuintegralmapnetwork}
Define the $L \times L$ matrix $M$ via
\begin{align*}
M_{kl} &=
\begin{cases}
\mu_{lk}^{\pi} & \text{ if } k \in \left\lbrace 1 , \dotsc , L \right\rbrace, l \in \mathcal{I}_{k}; \\
-\sum_{l \in \mathcal{I}_{k}} \mu_{kl}^{\pi} & \text{ if } k \in \left\lbrace 1 , \dotsc , L \right\rbrace, l = k.
\end{cases}
\end{align*}
Let $h \colon \mathbb{R}^{L} \to \mathbb{R}^{L}$ be given by $h \left( a \right) = Ma$. Then $h$ is Lipschitz continuous and Eq.\ \eqref{eq:ctuintegralmap} defines a function $H \colon \mathbb{R}^{L} \times D \left( \left[ 0,\infty \right) ; \mathbb{R}^{L} \right) \to D \left( \left[ 0,\infty \right) ; \mathbb{R}^{L} \right)$ that has the continuity properties described in Lemma~\ref{lem:ctuintegralmap}. If $y = H \left( b,x \right)$, then
\begin{align}
y_{k} \left( t \right)
&= b_{k} + x_{k} \left( t \right) + \sum_{l \in \mathcal{I}_{k}} \int_{0}^{t} \mu_{lk}^{\pi} y_{l} \left( s \right) \, \dd s - \sum_{l \in \mathcal{I}_{k}} \int_{0}^{t} \mu_{kl}^{\pi} y_{k} \left( s \right) \, \dd s. \label{eq:ctuintegralmapnetwork}
\end{align}
\end{corollary}

Throughout this section, the following conditions are in force,  with
\begin{align}
\beta = \max \left\lbrace 1/2 , 1 - \alpha/2 \right\rbrace.
\label{eq:betadef}
\end{align}

\begin{assumption}\label{AS2}
In the $n$-th system,
\begin{itemize}
\item[$\circ$]
with $Q^{n} \left( 0 \right)$ denoting the vector with the initial number of jobs in queue $k$, we assume  that 
\begin{align}
n^{1 - \beta} \left( \frac{1}{n} Q^{n} \left( 0 \right) - \rho \left( 0 \right) \right) \to X \left( 0 \right) \label{eq:initcondconv}
\end{align}
in distribution, where $X \left( 0 \right)$ is some random variable;
\item[$\circ$] with $\lambda_{k}^{n}$ the arrival rate to queue $k$, we assume that
\begin{align}
\lim_{n \to \infty} n^{1 - \beta} \left( \frac{1}{n} \lambda^{n}_{k} - \lambda_{k} \right) = \hat{\lambda}_{k};\label{eq:lambdahatdef}
\end{align}
\item[$\circ$] with $\mu_{kl}^{n}$
 the service speed for jobs to be routed from queue $k$ to queue $l$, we assume that
 \begin{align}
\lim_{n \to \infty} n^{1 - \beta} \left( \mu^{n}_{kl} - \mu_{kl} \right) = \hat{\mu}_{kl}; \label{eq:muhatdef}
\end{align}
\item[$\circ$] 
the background process, denoted by $J^{n}$, corresponds to a continuous-time Markov chain having irreducible generator matrix $n^{\alpha} \genQ$ (for some $\alpha>0$) and invariant distribution $\pi$.
\end{itemize}
\end{assumption}

Note that the assumption concerning the initial number of jobs is compatible with $\tfrac{1}{n} Q^{n} \left( 0 \right)$ converging in probability to $\rho \left( 0 \right)$, due to Slutsky's Lemma and $\beta - 1$ being strictly negative. In fact, for these reasons \eqref{eq:initcondconv} implies that $\tfrac{1}{n} Q^{n} \left( 0 \right)$ converges in probability to $\rho \left( 0 \right)$.

The process of interest here is an an appropriately centered and normalized version of $Q^n$, namely the stochastic process $\hat{Q}^{n}$, which is defined via
\begin{align}
\hat{Q}_{k}^{n} \left( t \right) &= n^{1 - \beta} \left( \frac{1}{n} Q_{k}^{n} \left( t \right) - \rho_{k} \left( t \right) \right). \label{eq:Qhatdeforig}
\end{align}
In words, $\hat{Q}^{n}$ is the scaled and centered $\mathbb{R}^{L}$-valued process describing how the number of jobs in the network fluctuates around the fluid limit. As announced before, our objective is to show that $\hat{Q}^{n}$ converges to a diffusion process. Anticipating the application of continuous-mapping methods, we represent $\hat{Q}^{n}$ as
\begin{align}
\begin{split}
\hat{Q}_{k}^{n} \left( t \right) 
&= \hat{Q}_{k}^{n} \left( 0 \right) + \hat{X}^{n}_{k} \left( t \right) + \sum_{l \in \mathcal{I}_{k}} \int_{0}^{t} \mu_{lk} \left( J^{n} \left( s \right) \right) \hat{Q}_{l}^{n} \left( s \right) \, \dd s - \sum_{l \in \mathcal{I}_{k}} \int_{0}^{t} \mu_{kl} \left( J^{n} \left( s \right) \right) \hat{Q}_{k}^{n} \left( s \right) \, \dd s,
\end{split} \label{eq:Qhatdef}
\end{align}
where
\begin{align}
\hat{X}^{n}_{k} \left( t \right) &= n^{1-\beta} \bar{X}^{n}_{k} \left( t \right)
\label{eq:Xhatdef}
\end{align}
and
\begin{align}
\begin{split}
\bar{X}^{n}_{k} \left( t \right) &= \bar{X}^{n}_{1,k} \left( t \right) + \bar{X}^{n}_{2,k} \left( t \right) + \bar{X}^{n}_{3,k} \left( t \right) + \sum_{l \in \mathcal{I}_{k}} \left( \bar{X}^{n}_{4,lk} \left( t \right) + \bar{X}^{n}_{5,lk} \left( t \right) + \bar{X}^{n}_{6,lk} \left( t \right) \right) \\
&\phantom{{} =} {} - \sum_{l \in \mathcal{I}_{k}} \left( \bar{X}^{n}_{4,kl} \left( t \right) + \bar{X}^{n}_{5,kl} \left( t \right) + \bar{X}^{n}_{6,kl} \left( t \right) \right),
\end{split} \label{eq:Xhatinputdef}
\end{align}
with the processes $\bar{X}^{n}_{1,k}$, $\bar{X}^{n}_{2,k}$, $\bar{X}^{n}_{3,k}$, $\bar{X}^{n}_{4,kl}$, $\bar{X}^{n}_{5,kl}$, and $\bar{X}^{n}_{6,kl}$ as defined in Eq.\ \eqref{eq:fluidproofXbardef} (see Section \ref{sec:fluidproof}).

Observe that the term $\mu_{kl} \left( J^{n} \left( s \right) \right)$ in the integral $\int_{0}^{t} \mu_{lk} \left( J^{n} \left( s \right) \right) \hat{Q}_{l}^{n} \left( s \right) \, \dd s$ above makes a direct application of Lemma~\ref{lem:ctuintegralmap} impossible. We deal with this problem by first showing that $\hat{Q}^{n}$ is asymptotically equivalent to the stochastic process $\tilde{Q}^{n}$, which is defined via
\begin{align}
\tilde{Q}_{k}^{n} \left( t \right) 
&= \hat{Q}_{k}^{n} \left( 0 \right) + \hat{X}^{n}_{k} \left( t \right) + \sum_{l \in \mathcal{I}_{k}} \int_{0}^{t} \mu_{lk}^{\pi} \tilde{Q}_{l}^{n} \left( s \right) \, \dd s - \sum_{l \in \mathcal{I}_{k}} \int_{0}^{t} \mu_{kl}^{\pi} \tilde{Q}_{k}^{n} \left( s \right) \, \dd s
\label{eq:Qtildekdef}
\end{align}
or, equivalently,
\begin{align}
\tilde{Q}^{n}
&= H \left( \hat{Q}^{n} \left( 0 \right) , \hat{X}^{n} \right),
\label{eq:QtildekdefviaH}
\end{align}
with $H$ denoting the continuous integral map from Example \ref{ex:ctuintegralmapnetwork}.
Observe that (\ref{eq:Qtildekdef})  differs from  (\ref{eq:Xhatinputdef})  in that the modulated server speeds $\mu_{kl} \left( J^{n} \left( s \right) \right)$ are replaced by their time-average counterparts.
The equivalence of both systems is motivated by the following intuitive reasoning.

Our idea is to establish the diffusion result relying on classical continuous-mapping arguments: we would like to view $\hat{Q}^{n}$ as a continuous integral map of some input processes. If the server speeds were not modulated, it would be clear that the input is given by $\hat{Q}^{n} \left( 0 \right)$ and $\hat{X}^{n}$, and that $\hat{Q}^{n}$ is a continuous integral map of these input processes. Then everything would fall nicely within the scope of Lemma~\ref{lem:ctuintegralmap}. Unfortunately, the server speeds do depend on the background process, so this approach does not work.

However, notice that the service requirements are not scaled, whereas the background process is scaled proportionally to $n^\alpha$. Intuitively, this means that a service requirement remains fixed as $n$ becomes large, while the background Markov chain is jumping faster and faster and approaches its equilibrium distribution. Thus, as $n$ becomes large, a service requirement is expected to effectively experience the average or equilibrium server speed.
These arguments suggest that, for large $n$, we may replace the modulated server speeds $\mu_{kl} \left( J^{n} \left( s \right) \right)$ by the average server speeds $\mu_{kl}^{\pi}$. More precisely, this means that the processes $\hat{Q}^{n}$ and $\tilde{Q}^{n}$ should be asymptotically equivalent.

There are two major flaws in the reasoning above. First, $\hat{Q}^{n}$ is a scaled and centered queueing process rather than an actual queueing process; it is unclear what a `service requirement' is in this context. Second, and more importantly, the term $\mu_{kl} \left( J^{n} \left( s \right) \right)$ also appears in the input via the $\hat{X}^{n}$ (see Eq.\ \eqref{eq:fluidproofXbardef} in Section \ref{sec:fluidproof}). So why not take the average server speed there? The not so convincing answer is that in the input $\mu_{kl} \left( J^{n} \left( s \right) \right)$ does not act as a server speed but is more like a measure of how the parameters deviate from their averaged counterparts.

Again, the intuitive reasoning above is not entirely convincing from a mathematical point of view, but at least it has given us a potentially useful idea. We formalize that idea in the following lemma.

\begin{lemma}
\label{lem:asymptoticequivalence}
The process $\tilde{Q}^{n}$ is asymptotically equivalent to $\hat{Q}^{n}$, meaning that $\hat{Q}^{n} - \tilde{Q}^{n}$ converges uoc in probability to the zero process as $n \to \infty$.
\end{lemma}

This lemma is useful because it has the following important corollary as an immediate consequence. Actually, the next result is the very reason for introducing $\tilde{Q}^{n}$: it claims that proving weak convergence of $\hat{Q}^{n}$ is equivalent to proving weak convergence of the simpler process $\tilde{Q}^{n}$, and that both processes must have the same limit.

\begin{corollary}
The process $\hat{Q}^{n}$ converges weakly to $X$ if and only if $\tilde{Q}^{n}$ converges weakly to $X$.
\end{corollary}

Thus, to prove weak convergence of $\hat{Q}^{n}$ it suffices to prove weak convergence of $\tilde{Q}^{n}$. However, because the process $\tilde{Q}^{n}$ was defined as a continuous function of some input process (cf.\ Eq.\ \eqref{eq:QtildekdefviaH}), we only have to prove weak convergence of this input process. These observations are the basis of the proof of the following theorem, which is the main result of this section.

The statement of the theorem contains two equivalent systems of stochastic integral equations. In the first system, each term corresponds to the term in Eq.\ \eqref{eq:Qhatdef} or Eq.\ \eqref{eq:Xhatinputdef} from which is is derived. In the second system, we have reordered the terms for easier interpretation and to show how the form of the integral map $H$ reappears in different parts of the limit.

\begin{theorem}
\label{thm:diffusionlimit}
Under Assumption \ref{AS2}, the processes $\tilde{Q}^{n}$ and $\hat{Q}^{n}$ both converge weakly to the unique solution $X$ of the system of stochastic integral equations
\begin{align}
\begin{split}
X_{k} \left( t \right)
&= X_{k} \left( 0 \right) + \ind{\alpha \geq 1} \sqrt{\lambda_{k}^{\pi}} B_{1,k} \left( t \right) + \hat{\lambda}_{k}^{\pi} t + \ind{\alpha \leq 1} \lambda_{k}^{\transpose} B_{\mathrm{b}} \left( t \right) \\
&\phantom{{} = {}} {} + \sum_{l \in \mathcal{I}_{k}} \Bigg( \ind{\alpha \geq 1} \int_{0}^{t} \sqrt{ \mu_{lk}^{\pi} \rho_{l} \left( s \right) } \, \dd B_{lk} \left( s \right) + \int_{0}^{t} \hat{\mu}_{lk}^{\pi} \rho_{l} \left( s \right) \, \dd s + \ind{\alpha \leq 1} \mu_{lk}^{\transpose} \int_{0}^{t} \rho_{l} \left( s \right) \, \dd B_{\mathrm{b}} \left( s \right) \Bigg) \\
&\phantom{{} = {}} {} - \sum_{l \in \mathcal{I}_{k}} \Bigg( \ind{\alpha \geq 1} \int_{0}^{t} \sqrt{ \mu_{kl}^{\pi} \rho_{k} \left( s \right) } \, \dd B_{kl} \left( s \right) + \int_{0}^{t} \hat{\mu}_{kl}^{\pi} \rho_{k} \left( s \right) \, \dd s + \ind{\alpha \leq 1} \mu_{kl}^{\transpose} \int_{0}^{t} \rho_{k} \left( s \right) \, \dd B_{\mathrm{b}} \left( s \right) \Bigg) \\
&\phantom{{} = {}} {} + \sum_{l \in \mathcal{I}_{k}} \int_{0}^{t} \mu_{lk}^{\pi} X_{l} \left( s \right) \dd s - \sum_{l \in \mathcal{I}_{k}} \int_{0}^{t} \mu_{kl}^{\pi} X_{k} \left( s \right) \dd s
\end{split}
\end{align}
or, equivalently,
\begin{align}
\begin{split}
X_{k} \left( t \right) 
&= X_{k} \left( 0 \right) + \int_{0}^{t} \left( \hat{\lambda}_{k}^{\pi} + \sum_{l \in \mathcal{I}_{k}} \hat{\mu}_{lk}^{\pi} \rho_{l} \left( s \right) - \sum_{l \in \mathcal{I}_{k}} \hat{\mu}_{kl}^{\pi} \rho_{k} \left( s \right) \right) \, \dd s \\
&\phantom{{} = {}} {} + \ind{\alpha \leq 1} \int_{0}^{t} \left( \lambda_{k}^{\transpose} + \sum_{l \in \mathcal{I}_{k}} \mu_{lk}^{\transpose} \rho_{l} \left( s \right) - \sum_{l \in \mathcal{I}_{k}} \mu_{kl}^{\transpose} \rho_{k} \left( s \right) \right) \, \dd B_{\mathrm{b}} \left( s \right) \\
&\phantom{{} = {}} {} + \ind{\alpha \geq 1} \int_{0}^{t} \Bigg( \sqrt{ \lambda_{k}^{\pi} } \, \dd B_{1,k} \left( s \right) + \sum_{l \in \mathcal{I}_{k}} \sqrt{ \mu_{lk}^{\pi} \rho_{l} \left( s \right) } \, \dd B_{lk} \left( s \right) - \sum_{l \in \mathcal{I}_{k}} \sqrt{ \mu_{kl}^{\pi} \rho_{k} \left( s \right) } \, \dd B_{kl} \left( s \right) \Bigg) \\
&\phantom{{} = {}} {} + \sum_{l \in \mathcal{I}_{k}} \int_{0}^{t} \mu_{lk}^{\pi} X_{l} \left( s \right) \dd s - \sum_{l \in \mathcal{I}_{k}} \int_{0}^{t} \mu_{kl}^{\pi} X_{k} \left( s \right) \dd s.
\end{split}
\end{align}
Here, the processes $B_{1,k}$ and $B_{kl}$ are independent standard Brownian motions. The process $B_{\mathrm{b}}$ is an independent $d$-dimensional Brownian motion with $\prequadvar{B_{\mathrm{b}}} \left( t \right) = \Sigma  t$, where $\Sigma:=\diag{\pi} D + D^{\transpose} \diag{\pi}$ and the matrix $D$ denotes the deviation matrix corresponding to the transition rate matrix $\genQ$ (as defined in Section \ref{sec:stateoccmeas}).
\end{theorem}


\section{Proofs and auxiliary results}\label{proofs}
In this section, we provide the proofs of the three main results, namely the fluid limit (Lemma \ref{lem:fluidlimit}), asymptotic equivalence (Lemma \ref{lem:asymptoticequivalence}), and the diffusion limit (Theorem \ref{thm:diffusionlimit}). Additionally, we present several auxiliary results which are used in those proofs. This section is organized in the following way.

First, we summarize some basic properties concerning the convergence of the state occupation measure of a Markov chain. These properties are essential ingredients of the upcoming proofs. We start with the proof of the fluid limit. This result is derived under slightly weaker conditions than the diffusion limit. Moreover, it is used in the other two proofs. Then, we proceed by establishing the asymptotic equivalence result in Lemma \ref{lem:asymptoticequivalence}. Finally, relying on both the fluid limit and the asymptotic equivalence result, we give the proof of the diffusion limit.

\subsection{State occupation measures of Markov chains}
\label{sec:stateoccmeas}
In the three upcoming proofs, the behavior of the state occupation measure of a scaled Markov chain plays a key role. In particular, we exploit weak convergence results for stochastic integrals with respect to the state occupation measure. To make this paper more or less self-contained, we collect the most important results concerning state occupation measures here. Their proofs can be found in \cite{Jansen2017X}.

As usual, $J$ denotes a continuous-time Markov chain with state space $\left\lbrace 1 , \dotsc , d \right\rbrace$ for some $d \in \mathbb{N}$. It has $d \times d$ generator matrix $\genQ$, which we assume to be irreducible. We denote the stationary distribution corresponding to $\genQ$ by $\pi$, i.e., $\pi$ is the unique $d \times 1$ vector (with nonnegative entries summing to $1$) that solves the equation $\pi^{\transpose} \genQ = 0$. Additionally, $D$ denotes the deviation matrix corresponding to $\genQ$. The entries of $D$ are given by
\begin{align*}
D_{ij} = \int_{0}^{\infty} \left( \mathbb{P} \left( J \left( s \right) = j \, \middle\vert \, J \left( 0 \right) = i \right) - \pi_{j} \right) \, \dd s.
\end{align*}
For some background on deviation matrices, see \cite{CV2002}. Define $\Sigma= \diag{\pi} D + D^{\transpose} \diag{\pi}$.

The state indicator function of $J$ is the $\mathbb{R}^{d}$-valued function $K$ defined via
\begin{align*}
K \left( i ; t \right)
&= \ind{J \left( t \right) = i}
\end{align*}
for $i \in \left\lbrace 1 , \dotsc , d \right\rbrace$ and $t \geq 0$. Then the state occupation measure corresponding to $J$ is given by
\begin{align*}
\int_{0}^{t} K \left( s \right) \, \dd s.
\end{align*}
Thus, the state indicator function $K$ is used to keep track of the state of $J$ at a given time, whereas the state occupation measure records how much time $J$ has spent in each state during a time interval.

The following theorem concerns the behavior of the state occupation measure of a Markov chain when the generator matrix is scaled. Let $\genQ$,  $\pi$, $D$, and $\Sigma$ be as defined above. Let $J^{n}$ be a continuous-time Markov chain with generator matrix $n^{\alpha} \genQ$ (for some $\alpha>0$) and state indicator function $K^{n}$. Define the stochastic process $G^{n}$ via
\begin{align}
G^{n} \left( t \right) = \int_{0}^{t} \left( K^{n} \left( s \right) - \pi \right) \, \dd s.
\label{eq:defGn}
\end{align}
Additionally, let $Y^{n}$ denote the Dynkin martingale corresponding to $K^{n}$ (cf.\ \cite[Lem.~2.6.18]{ae2004}).

\begin{theorem}
\label{thm:occupationmeasureconvergence}
The process $G^{n}$ converges to the zero process uoc in probability for $n \to \infty$. Additionally, $n^{\alpha / 2} G^{n}$ converges weakly to a Brownian motion $B_{\mathrm{b}}$ having predictable quadratic variation process $\prequadvar{B_{\mathrm{b}}} \left( t \right) = \Sigma t$.
\end{theorem}

In the proof of the diffusion limit, we also need weak convergence results of the form
\begin{align*}
H^{n} \stdot n^{\alpha / 2} G^{n} \Rightarrow H \stdot G,
\end{align*}
where $Y \stdot X$ denotes the It\^{o} integral of $Y$ with respect to $X$. This type of convergence generally does not hold, even when $G$ is a Brownian motion, $H$ is the zero process, and $H^{n}$ converges to $H$ uoc almost surely. However, under even stronger conditions this type of convergence does hold. We state such conditions in the next theorem.

Before we can state the theorem, we need to introduce some more notation (which by and large follows \cite[p.\ 204]{js2003}). Let $X$ be a $d$-dimensional locally square-integrable martingale with respect to a filtration $\mathds{F}$. For simplicity, we assume that
\begin{align*}
\prequadvar{X} \left( t \right)
&= \int_{0}^{t} C \left( s \right) \, \dd s
\end{align*}
for a predictable process $C$ taking values in the set of all symmetric nonnegative $d \times d$ matrices. We denote by $L^{2}_{\mathrm{loc}} \left( X \right)$ the set of predictable processes $H$ taking values in $\mathbb{R}^{k \times d}$ such that the process
\begin{align*}
\int_{0}^{t} H \left( s \right) C \left( s \right) H \left( s \right)^{\transpose} \, \dd s
\end{align*}
is locally integrable.

\begin{theorem}
\label{thm:markovintegralconvergence}
Consider the setting  of Theorem \ref{thm:occupationmeasureconvergence}. For fixed $m \in \mathbb{N}$, let $H^{1,n} , \dotsc , H^{m,n}$ be c\`{a}dl\`{a}g processes in $L^{2}_{\mathrm{loc}} \left( Y^{n} \right)$ and define $H^{j,n}_{-}$ via $H^{j,n}_{-} \left( t \right) = H^{j,n} \left( t- \right)$. 
Assume that each $H^{j,n}$ converges to a deterministic continuous function $H^{j}$ uoc in probability and that $n^{- \alpha / 2} H^{j,n}$ is a finite variation process whose total variation process converges to the zero process uoc in probability. Then $\big( H^{1,n}_{-} \stdot n^{\alpha / 2} G^{n} , \dotsc , H^{m,n}_{-} \stdot n^{\alpha / 2} G^{n} \big)$ converges weakly to $\big( H^{1} \stdot B_{\mathrm{b}} , \dotsc , H^{m} \stdot B_{\mathrm{b}} \big)$.
\end{theorem}

\subsection{Proof of the fluid limit}
\label{sec:fluidproof}
To establish the fluid limit in Lemma \ref{lem:fluidlimit}, we follow a standard approach by combining the Functional Law of Large Numbers (FLLN) for Poisson processes with an application of Gronwall's Lemma. The fact that the parameters are modulated by a background process is a minor complication here, but we can deal with that by exploiting convergence properties of irreducible Markov chains (cf.\ Section \ref{sec:stateoccmeas}).

Impose Assumption \ref{AS1}, as introduced in Section \ref{sec:fluidlimit}. Under these conditions, Lemma \ref{lem:fluidlimit} claims uoc convergence of $\bar{Q}_{k}^{n}$ to the zero process, where
\begin{align*}
\bar{Q}_{k}^{n} \left( t \right) &= \frac{1}{n} Q_{k}^{n} \left( t \right) - \rho_{k} \left( t \right).
\end{align*}
We prove this claim in the following way.

Using the expressions in Eq.~\eqref{eq:scaledsysdef} and Eq.~\eqref{eq:rhodef}, some simple algebra gives us
\begin{align*}
\bar{Q}_{k}^{n} \left( t \right)
&= \bar{Q}_{k}^{n} \left( 0 \right) + \bar{X}^{n}_{1,k} \left( t \right) + \bar{X}^{n}_{2,k} \left( t \right) + \bar{X}^{n}_{3,k} \left( t \right) \\
&\phantom{{} = {}} {} + \sum_{l \in \mathcal{I}_{k}} \left( \bar{X}^{n}_{4,lk} \left( t \right) + \bar{X}^{n}_{5,lk} \left( t \right) + \int_{0}^{t} \mu_{lk} \left( J^{n} \left( s \right) \right) \bar{Q}_{l}^{n} \left( s \right) \, \dd s + \bar{X}^{n}_{6,lk} \left( t \right) \right) \\
&\phantom{{} = {}} {} - \sum_{l \in \mathcal{I}_{k}} \left( \bar{X}^{n}_{4,kl} \left( t \right) + \bar{X}^{n}_{5,kl} \left( t \right) + \int_{0}^{t} \mu_{kl} \left( J^{n} \left( s \right) \right) \bar{Q}_{k}^{n} \left( s \right) \, \dd s + \bar{X}^{n}_{6,kl} \left( t \right) \right),
\end{align*}
with
\begin{align}
\begin{split}
\bar{X}^{n}_{1,k} \left( t \right) &= \frac{1}{n} A_{k} \left( n \int_{0}^{t} \frac{1}{n} \lambda_{k}^{n} \left( J^{n} \left( s \right) \right) \, \dd s \right) - \int_{0}^{t} \frac{1}{n} \lambda_{k}^{n} \left( J^{n} \left( s \right) \right) \, \dd s, \\
\bar{X}^{n}_{2,k} \left( t \right) &= \int_{0}^{t} \frac{1}{n} \lambda_{k}^{n} \left( J^{n} \left( s \right) \right) \, \dd s - \int_{0}^{t} \lambda_{k} \left( J^{n} \left( s \right) \right) \, \dd s, \\
\bar{X}^{n}_{3,k} \left( t \right) &= \int_{0}^{t} \lambda_{k} \left( J^{n} \left( s \right) \right) \, \dd s - \int_{0}^{t} \lambda_{k}^{\pi} \, \dd s, \\
\bar{X}^{n}_{4,kl} \left( t \right) &= \frac{1}{n} S_{kl} \left( n \int_{0}^{t} \mu_{kl}^{n} \left( J^{n} \left( s \right) \right) \frac{1}{n} Q_{k}^{n} \left( s \right) \, \dd s \right) - \int_{0}^{t} \mu_{kl}^{n} \left( J^{n} \left( s \right) \right) \frac{1}{n} Q_{k}^{n} \left( s \right) \, \dd s, \\
\bar{X}^{n}_{5,kl} \left( t \right) &= \int_{0}^{t} \left( \mu_{kl}^{n} \left( J^{n} \left( s \right) \right) - \mu_{kl} \left( J^{n} \left( s \right) \right) \right) \frac{1}{n} Q_{k}^{n} \left( s \right) \, \dd s, \\
\bar{X}^{n}_{6,kl} \left( t \right) &= \int_{0}^{t} \left( \mu_{kl} \left( J^{n} \left( s \right) \right) - \mu_{kl}^{\pi} \right) \rho_{k} \left( s \right) \, \dd s.
\end{split}
\label{eq:fluidproofXbardef}
\end{align}
Now suppose that all six types of the above processes converge to the zero process uoc in probability when $\tfrac{1}{n} Q^{n} \left( 0 \right)$ converges to $\rho \left( 0 \right)$ in probability. Then, for fixed $T \geq 0$, the random variables
\begin{align*}
\bar{Y}^{n}_{1} \left( T \right) &= \sum_{k=1}^{L} \left( \sup_{t \in \left[ 0,T \right]} \abs{ \bar{X}^{n}_{1,k} \left( t \right) } + \sup_{t \in \left[ 0,T \right]} \abs{ \bar{X}^{n}_{2,k} \left( t \right) } + \sup_{t \in \left[ 0,T \right]} \abs{ \bar{X}^{n}_{3,k} \left( t \right) } \right), \\
\bar{Y}^{n}_{2} \left( T \right) &= \sum_{k=1}^{L} \sum_{l \in \mathcal{I}_{k}} \left( \sup_{t \in \left[ 0,T \right]} \abs{ \bar{X}^{n}_{4,lk} \left( t \right) } + \sup_{t \in \left[ 0,T \right]} \abs{ \bar{X}^{n}_{5,lk} \left( t \right) } + \sup_{t \in \left[ 0,T \right]} \abs{ \bar{X}^{n}_{6,lk} \left( t \right) } \right), \\
\intertext{and}
\bar{Y}^{n}_{3} \left( T \right) &= \sum_{k=1}^{L} \sum_{l \in \mathcal{I}_{k}} \left( \sup_{t \in \left[ 0,T \right]} \abs{ \bar{X}^{n}_{4,kl} \left( t \right) } + \sup_{t \in \left[ 0,T \right]} \abs{ \bar{X}^{n}_{5,kl} \left( t \right) } + \sup_{t \in \left[ 0,T \right]} \abs{ \bar{X}^{n}_{6,kl} \left( t \right) } \right)
\end{align*}
converge to $0$ in the same way. Observe that
\begin{align*}
\sum_{k=1}^{L} \abs{ \bar{Q}_{k}^{n} \left( t \right) }
&\leq \sum_{k=1}^{L} \abs{ \bar{Q}_{k}^{n} \left( 0 \right) } + \bar{Y}^{n}_{1} \left( T \right) + \bar{Y}^{n}_{2} \left( T \right) + \bar{Y}^{n}_{3} \left( T \right) \\
&\phantom{{} \leq{} } {} + \sum_{k=1}^{L} \sum_{l \in \mathcal{I}_{k}} \int_{0}^{t} \mu^{*} \abs{ \bar{Q}_{l}^{n} \left( s \right) } \, \dd s + \sum_{k=1}^{L} \sum_{l \in \mathcal{I}_{k}} \int_{0}^{t} \mu^{*} \abs{ \bar{Q}_{k}^{n} \left( s \right) } \, \dd s
\end{align*}
for all $t \in \left[ 0,T \right]$, where $\mu^{*} = \sum_{k=1}^{L} \sum_{l \in \mathcal{I}_{k}} \sum_{i=1}^{d} \mu_{kl} \left( i \right)$. It follows that
\begin{align*}
\sum_{k=1}^{L} \abs{ \bar{Q}_{k}^{n} \left( t \right) }
&\leq \sum_{k=1}^{L} \abs{ \bar{Q}_{k}^{n} \left( 0 \right) } + \bar{Y}^{n}_{1} \left( T \right) + \bar{Y}^{n}_{2} \left( T \right) + \bar{Y}^{n}_{3} \left( T \right) + \int_{0}^{t} 2 \left( L-1 \right) \mu^{*} \sum_{k=1}^{L} \abs{ \bar{Q}_{k}^{n} \left( s \right) } \, \dd s
\end{align*}
for all $t \in \left[ 0,T \right]$. An application of Gronwall's Lemma (cf.\ \cite[Th.~A.5.1]{ek1986}) then leads to
\begin{align*}
\sum_{k=1}^{L} \abs{ \bar{Q}_{k}^{n} \left( t \right) }
&\leq \left( \sum_{k=1}^{L} \abs{ \bar{Q}_{k}^{n} \left( 0 \right) } + \bar{Y}^{n}_{1} \left( T \right) + \bar{Y}^{n}_{2} \left( T \right) + \bar{Y}^{n}_{3} \left( T \right) \right) e^{2 \left( L-1 \right) \mu^{*} T}
\end{align*}
for all $t \in \left[ 0,T \right]$, so $\sup_{t \in \left[ 0,T \right]} \sum_{k=1}^{L} \abs{ \bar{Q}_{k}^{n} \left( t \right) }$ converges to $0$ in probability under the assumptions imposed (in particular, recall that $\bar Q^n(0)\to 0$ in
probability).

Consequently, it remains to show that $\bar{Y}^{n}_{1} \left( T \right)$, $\bar{Y}^{n}_{2} \left( T \right)$, and $\bar{Y}^{n}_{3} \left( T \right)$ converge to $0$ in probability as soon as $\tfrac{1}{n} Q^{n} \left( 0 \right)$ converges to $\rho \left( 0 \right)$ in probability.
To this end, 
we first define \[\lambda^{*} = 1 + \sum_{k=1}^{L} \sum_{i=1}^{d} \lambda_{k} \left( i \right).\] Then $\int_{0}^{t} \tfrac{1}{n} \lambda_{k}^{n} \left( J^{n} \left( s \right) \right) \, \dd s \leq 2 \lambda^{*} t$ for all $n$ large enough. Moreover, $\tfrac{1}{n} A_{k} \left( n\, 2 \lambda^{*} t \right)$ converges uoc to $2 \lambda^{*} t$ almost surely by the FLLN for Poisson processes (cf.\ \cite[Th.\ 5.5.10]{cy2001}), implying that $\sup_{t \in \left[ 0,T \right]} \absj{ \bar{X}_{1,k}^{n} \left( t \right) }$ converges to $0$ almost surely.

Next, observe that $\sup_{t \in \left[ 0,T \right]} \absj{ \bar{X}_{2,k}^{n} \left( t \right) }$ converges to $0$ almost surely, due to the convergence of $\tfrac{1}{n} \lambda_{k}^{n}$ to $\lambda_{k}$. Also observe that both $\sup_{t \in \left[ 0,T \right]} \absj{ \bar{X}_{3,k}^{n} \left( t \right) }$ and $\sup_{t \in \left[ 0,T \right]} \absj{ \bar{X}_{6,kl}^{n} \left( t \right) }$ converge to $0$ in probability by Theorem \ref{thm:occupationmeasureconvergence} and Theorem \ref{thm:markovintegralconvergence}, respectively.

To analyze $\bar{X}_{4,kl}^{n} \left( t \right)$ and $\bar{X}_{5,kl}^{n} \left( t \right)$, we need to know more about the behavior of $Q_{k}^{n} \left( s \right)$. Clearly, the queueing dynamics are such that
\begin{align*}
\abs{ \frac{1}{n} Q_{k}^{n} \left( s \right) } \leq \sum_{m \in \mathcal{I}} \left( \frac{1}{n} Q_{m}^{n} \left( 0 \right) + \frac{1}{n} A_{m} \left( n \lambda^{*} T \right) \right)
\end{align*}
for all $s \in \left[ 0,T \right]$. The sum on the right-hand side converges to $\sum_{m \in \mathcal{I}} \left( \rho_{m} \left( 0 \right) + \lambda^{*} T \right)$ in probability, due to Slutsky's Lemma (cf.\ \cite[Lem.\ 2.8]{vandervaart1998}). Another application of Slutsky's Lemma shows that 
\begin{align*}
\abs{ \mu_{kl}^{n} \left( J^{n} \left( s \right) \right) - \mu_{kl} \left( J^{n} \left( s \right) \right) } \sum_{l=1}^{L} \left( \frac{1}{n} Q_{l}^{n} \left( 0 \right) + \frac{1}{n} A_{l} \left( n \lambda^{*} T \right) \right)
\end{align*}
converges to $0$ in probability. Then $\sup_{t \in \left[ 0,T \right]} \absj{ \bar{X}_{5,kl}^{n} \left( t \right) }$ must converge to $0$ in probability, too. Regarding $\bar{X}_{4,kl}^{n} \left( t \right)$, we observe that, for fixed $\eps > 0$, we have the very crude inequality
\begin{align*}
\int_{0}^{t} \mu_{kl}^{n} \left( J^{n} \left( s \right) \right) \frac{1}{n} Q_{k}^{n} \left( s \right) \, \dd s 
\leq \mu^{*} \sum_{m \in \mathcal{I}} \left( \rho_{m} \left( 0 \right) + \lambda^{*} T \right) T + \eps
\end{align*}
on a set with probability at least $1 - \eps$ for all $n$ large enough. Combined with the FLLN for the standard Poisson process $S_{kl}$, this gives us the convergence of $\sup_{t \in \left[ 0,T \right]} \absj{ \bar{X}_{4,kl}^{n} \left( t \right) }$ to $0$ in probability.

Using Slutsky's Lemma once again, we conclude that $\bar{Y}^{n}_{1} \left( T \right)$, $\bar{Y}^{n}_{2} \left( T \right)$, and $\bar{Y}^{n}_{3} \left( T \right)$ converge to $0$ in probability. This completes the proof.

\subsection{Proof of asymptotic equivalence}
\label{sec:asymptoticequivalenceproof}
Impose Assumption \ref{AS2}, as stated in Section \ref{sec:diffusionlimit}. Our goal is to prove the claim in Lemma \ref{lem:asymptoticequivalence}, namely that $\tilde{Q}^{n}$ is asymptotically equivalent to $\hat{Q}^{n}$. Recall that $\hat{Q}^{n}$ represents the scaled and centered queue content process and that $\tilde{Q}^{n}$ is a similar but simpler stochastic process defined using a continuous map. (Note that, in general, $\tilde{Q}^{n}$ is not a scaled and centered queue content process.)

For completeness, we first show that $\tilde{Q}^{n}$ is a well-defined stochastic process. 
We already know that $\hat{Q}^{n} \left( 0 \right)$ and $\hat{X}^{n}$ are well-defined random elements. As observed in Example \ref{ex:ctuintegralmapnetwork}, the integral map $H$ that defines $\tilde{Q}^{n}$ satisfies the conditions of Lemma~\ref{lem:ctuintegralmap}, so $\tilde{Q}^{n}$ is a continuous function (in the Skorokhod J$_{1}$ topology) of a stochastic process. Then $\tilde{Q}^{n}$ is a well-defined stochastic process.

We turn to proving asymptotic equivalence of $\hat{Q}^{n}$ and $\tilde{Q}^{n}$. The difference between these two processes can be written as
\begin{align*}
&\hat{Q}^{n}_{k} \left( t \right) - \tilde{Q}^{n}_{k} \left( t \right) \\
&\quad = \sum_{l \in \mathcal{I}_{k}} \left( \int_{0}^{t} \mu_{lk} \left( J^{n} \left( s \right) \right) \hat{Q}^{n}_{l} \left( s \right) \, \dd s - \int_{0}^{t} \mu_{lk}^{\pi} \tilde{Q}^{n}_{l} \left( s \right) \, \dd s \right) \\
&\phantom{\quad {} =} {} - \sum_{l \in \mathcal{I}_{k}} \left( \int_{0}^{t} \mu_{kl} \left( J^{n} \left( s \right) \right) \hat{Q}^{n}_{k} \left( s \right) \, \dd s - \int_{0}^{t} \mu_{kl}^{\pi} \tilde{Q}^{n}_{k} \left( s \right) \, \dd s \right) \\
&\quad = \sum_{l \in \mathcal{I}_{k}} \left( \int_{0}^{t} \left( \mu_{lk} \left( J^{n} \left( s \right) \right) - \mu_{lk}^{\pi} \right) \hat{Q}^{n}_{l} \left( s \right) \, \dd s  - \int_{0}^{t} \mu_{lk}^{\pi} \left[ \tilde{Q}^{n}_{l} \left( s \right) - \hat{Q}^{n}_{l} \left( s \right) \right] \, \dd s \right) \\
&\phantom{\quad {} =} {} - \sum_{l \in \mathcal{I}_{k}} \left( \int_{0}^{t} \left( \mu_{kl} \left( J^{n} \left( s \right) \right) - \mu_{kl}^{\pi} \right) \hat{Q}^{n}_{k} \left( s \right) \, \dd s  - \int_{0}^{t} \mu_{kl}^{\pi} \left[ \tilde{Q}^{n}_{k} \left( s \right) - \hat{Q}^{n}_{k} \left( s \right) \right] \, \dd s \right),
\end{align*}
so this difference is, in fact, given by the continuous integral map $H$ of some input process. Consequently, if this input process converges to the zero process uoc in probability, then the difference between $\hat{Q}^{n}$ and $\tilde{Q}^{n}$ must converge to the zero process uoc in probability, too.

This input process is essentially a sum of other input processes of the form
\begin{align}
\int_{0}^{t} \left( \mu_{kl} \left( J^{n} \left( s \right) \right) - \mu_{kl}^{\pi} \right) \hat{Q}^{n}_{k} \left( s \right) \, \dd s. \label{eq:MCintegral}
\end{align}
Thus, to prove the lemma, it suffices to establish the convergence of each of these terms to the zero process uoc in probability. 
We will do this in the remainder of this proof.

Recall from Section \ref{sec:stateoccmeas} that the state indicator function of $J^{n}$ is denoted by $K^{n}$ and that the Dynkin martingale associated with $J^{n}$ is denoted by $Y^{n}$.
Now rewrite the integral in Eq.~\eqref{eq:MCintegral} as
\begin{align}
\int_{0}^{t} n^{-\alpha/2} \hat{Q}^{n}_{k} \left( s \right) \mu_{kl}^{\transpose} n^{\alpha/2} \left( K^{n} \left( s \right) - \pi \right) \, \dd s
&= \int_{0}^{t} U^{n}_{kl} \left( s \right) n^{\alpha/2} \left( K^{n} \left( s \right) - \pi \right) \, \dd s, \label{eq:MCintegralrewrite}
\end{align}
where
\begin{align*}
U^{n}_{kl} \left( t \right)
&= n^{-\alpha/2} \hat{Q}^{n}_{k} \left( t- \right) \mu_{kl}^{\transpose}.
\end{align*}
Convergence of $U^{n}_{kl}$ to the zero process implies, under certain assumptions, convergence of the integral \eqref{eq:MCintegralrewrite} to the zero process; this is described Theorem \ref{thm:markovintegralconvergence}.

To be able to apply Theorem \ref{thm:markovintegralconvergence}, we have to check several properties. We will break this up in four steps and then apply Theorem \ref{thm:markovintegralconvergence} in the fifth and final step.

\emph{Step 1.}
We would like to verify (in this step and the next) that $U^{n}_{kl}$ is a process in $L^{2}_{\mathrm{loc}} \left( Y^{n} \right)$. First of all, we should have a filtration with respect to which $U^{n}_{kl}$ is predictable and $Y^{n}$ is a martingale. For this filtration, say $\mathds{G}$, we can just take the natural filtration of $J^{n}$ and add to this the entire filtration of all Poisson processes and $Q^{n} \left( 0 \right)$, which are independent of $J^{n}$. Then $Y^{n}$ is still a (locally square-integrable) martingale with respect to this filtration, due to the independence. Moreover, the left-continuous process $U^{n}_{kl}$ is adapted to $\mathds{G}$ and thus predictable. This settles the filtration and measurability issues of Theorem \ref{thm:markovintegralconvergence}.

\emph{Step 2.}
From the previous step, we know that $U^{n}_{kl}$ has the right measurability properties. To verify that $U^{n}_{kl}$ is a process in $L^{2}_{\mathrm{loc}} \left( Y^{n} \right)$, we also have to check that the process
\begin{align*}
\int_{0}^{t} U^{n}_{kl} \left( s \right) \left[ \diag{\genQ^{\transpose} K \left( s \right)} - \genQ^{\transpose} \diag{K \left( s \right)} - \diag{K \left( s \right)} \genQ \right] U^{n}_{kl} \left( s \right)^{\transpose} \, \dd s
\end{align*}
is locally integrable. This is rather obvious, as
\begin{align*}
\abs{ \hat{Q}^{n}_{k} \left( t \right) } 
&\leq n^{1 - \beta} \frac{1}{n} Q^{n}_{k} \left( t \right) + n^{1 - \beta} \rho_{k} \left( t \right) \\
&\leq n^{1 - \beta} \rho_{k} \left( t \right) + \sum_{l=1}^{L} A_{l} \left( \int_{0}^{t} \lambda_{l}^{n} \left( J^{n} \left( s \right) \right) \, \dd s \right) \\
&\leq n^{1 - \beta} \rho_{k} \left( t \right) + \sum_{l=1}^{L} A_{l} \left( \lambda^{*} t \right),
\end{align*}
where we now define $\lambda^{*}$ as $\sum_{k=1}^{L} \sum_{i=1}^{d} \lambda_{k} \left( i \right)$. From this, we conclude that $U^{n}_{kl}$ is a process in $L^{2}_{\mathrm{loc}} \left( Y^{n} \right)$.

\emph{Step 3.} Perhaps the most important step is to check that $U^{n}_{kl}$ converges to the zero process uoc in probability. From the definition of $\hat{Q}^{n}_{k}$ in Eq.\ \eqref{eq:Qhatdef}, it follows that
\begin{align}
U^{n}_{kl} \left( t \right) = n^{1 - \beta - \alpha/2} \left( \frac{1}{n} {Q}^{n}_{k} \left( t- \right) - \rho_{k} \left( t- \right) \right) \mu_{kl}^{\transpose}. \label{eq:Hnknonscale}
\end{align}
By the fluid limit in Lemma~\ref{lem:fluidlimit}, the process $\frac{1}{n} {Q}^{n}_{k} - \rho_{k}$ converges to the zero process uoc in probability. Because $1 - \beta - \alpha/2 = \min \left\lbrace 0 , \left( 1 - \alpha \right)/2 \right\rbrace$, we know that $U^{n}_{kl}$ converges to the zero process in the same way.

\emph{Step 4.}
Now we only have to verify that $n^{- \alpha/2} U^{n}_{kl}$ is a finite variation process whose total variation process converges to the zero process uoc in probability. From Eq.\ \eqref{eq:Hnknonscale}, we obtain
\begin{align}
n^{- \alpha/2} U^{n}_{kl} \left( t \right) 
&= n^{- \alpha/2} n^{1 - \beta - \alpha/2} \left( \frac{1}{n} {Q}^{n}_{k} \left( t- \right) - \rho_{k} \left( t- \right) \right) \mu_{kl}^{\transpose} \nonumber \\
&= n^{1 - \beta - \alpha} \left( \frac{1}{n} {Q}^{n}_{k} \left( t- \right) - \rho_{k} \left( t- \right) \right) \mu_{kl}^{\transpose}. \label{eq:Hnkscale}
\end{align}
Clearly, it suffices to show that $n^{1 - \beta - \alpha} \left( \frac{1}{n} {Q}^{n}_{k} \left( t- \right) - \rho_{k} \left( t- \right) \right)$ is a finite variation process whose total variation process converges to the zero process uoc in probability.

To this end, observe that $\rho_{k}$ is continuously differentiable, so it is of bounded variation. Moreover, $\rho_{k}$ does not depend on $n$ and we have 
\begin{align}
1 - \beta - \alpha = \min \left\lbrace - \alpha/2 , \left( 1 - 2\alpha \right)/2 \right\rbrace < 0, \label{eq:oneminusbetaminusalpha}
\end{align}
so the total variation process of $n^{1 - \beta - \alpha} \rho_{k}$ converges to the zero process uoc. Consequently, it suffices to show that $n^{1 - \beta - \alpha} \frac{1}{n} {Q}^{n}_{k} \left( t- \right)$ is a finite variation process whose total variation process converges to the zero process uoc in probability.

Now note that $\frac{1}{n} {Q}^{n}_{k} \left( t- \right)$ can be written as the difference of nondecreasing processes (cf.\ Eq.\ \eqref{eq:scaledsysdef}). This means that $n^{1 - \beta - \alpha} \frac{1}{n} {Q}^{n}_{k} \left( t- \right)$ is a finite variation process. It also follows from Eq.\ \eqref{eq:scaledsysdef} that the total variation process of $n^{1 - \beta - \alpha} \frac{1}{n} {Q}^{n}_{k} \left( t- \right)$ is given by
\begin{align*}
&n^{1 - \beta - \alpha} \frac{1}{n} Q_{k}^{n} \left( 0 \right) + n^{1 - \beta - \alpha} \frac{1}{n} A_{k} \left( n \int_{0}^{t} \frac{1}{n}\lambda_{k}^{n} \left( J^{n} \left( s \right) \right) \, \dd s \right) \\
&{} + n^{1 - \beta - \alpha} \sum_{l \in \mathcal{I}_{k}} \frac{1}{n} S_{lk} \left( n \int_{0}^{t} \mu_{lk}^{n} \left( J^{n} \left( s \right) \right) \frac{1}{n} Q_{l}^{n} \left( s \right) \, \dd s \right) \\
&{} + n^{1 - \beta - \alpha} \sum_{l \in \mathcal{I}_{k}} \frac{1}{n} S_{kl} \left( n \int_{0}^{t} \mu_{kl}^{n} \left( J^{n} \left( s \right) \right) \frac{1}{n} Q_{k}^{n} \left( s \right) \, \dd s \right).
\end{align*}
As a result of Lemma \ref{lem:fluidlimit} and the fact that $1 - \beta - \alpha = \min \left\lbrace - \alpha/2 , \left( 1 - 2\alpha \right)/2 \right\rbrace$, we obtain the convergence of this expression to the zero process uoc in probability.

As argued before, this implies that $n^{- \alpha/2} U^{n}_{kl}$ is a finite variation process whose total variation process converges to the zero process uoc in probability.

\emph{Step 5.}
Recall that we aim to show convergence of the integral in Eq.\ \eqref{eq:MCintegralrewrite} to the zero process via an application of Theorem \ref{thm:markovintegralconvergence}. In the previous four steps, we have verified the assumptions of Theorem \ref{thm:markovintegralconvergence}. Moreover, we have shown that $U^{n}_{kl}$ converges to the zero process uoc in probability, so Theorem \ref{thm:markovintegralconvergence} implies that the integral in Eq.\ \eqref{eq:MCintegralrewrite} converges to the zero process uoc in probability, as required.

\subsection{Proof of the diffusion limit}
\label{sec:diffusionproof}
Before we turn to the proof of the diffusion limit result of Theorem~\ref{thm:diffusionlimit}, let us first comment on notation. We will be mainly concerned with joint weak convergence of certain stochastic processes. Because this involves many different processes, we will not explicitly order them in a vector, so as to minimize potential confusion. As an example, we will just say that $Y_{k}$ and $Y_{kl}$ converge jointly, where it is tacitly understood that $k \in \left\lbrace 1 , \dotsc , L \right\rbrace$ and $l \in \mathcal{I}_{k}$ (the ranges of the indices should be clear from the context), and that all these processes are ordered in some vector.

Impose Assumption \ref{AS2}, as was stated in Section \ref{sec:diffusionlimit}. We would like to establish the diffusion limit in Theorem \ref{thm:diffusionlimit} under these conditions. To this end, it suffices to show that $\tilde{Q}^{n}$ converges weakly to the process $X$ given in Theorem \ref{thm:diffusionlimit}, due to Lemma \ref{lem:asymptoticequivalence}. However, because $\tilde{Q}^{n}$ is defined as a continuous function of certain input processes $\hat{Q}^{n} \left( 0 \right)$ and $\hat{X}^{n}$, the CMT implies that we only have to prove joint weak convergence of these input processes. Consequently, the remainder of this proof will focus on  establishing this joint weak convergence.

Following the notation from Section \ref{sec:stateoccmeas}, we denote by $K^{n}$ the state indicator function corresponding to the background process $J^{n}$. Recall that $K^{n}$ is a $\left\lbrace 0,1 \right\rbrace^{d}$-valued process that is given by $K^{n} \left( i ; t \right) = \ind{J^{n} \left( s \right) = i}$ for $i \in \left\lbrace 1 , \dotsc , d \right\rbrace$. 
In addition, we define the continuous semimartingale $G^{n}$ via
\begin{align*}
G^{n} \left( t \right)
&= \int_{0}^{t} \left( K^{n} \left( s \right) - \pi \right) \, \dd s.
\end{align*}

As argued before, we only have to prove weak convergence of the input processes of $\tilde{Q}^{n}$, i.e., we have to prove joint weak convergence of $\hat{Q}^{n} \left( 0 \right)$ and $\hat{X}^{n}$. We will show that weak convergence of $\hat{X}^{n}$ follows (after an application of the CMT) from the joint weak convergence of other, more basic processes. Joint convergence of $\hat{Q}^{n} \left( 0 \right)$ and these basic processes will be straightforward (as they are independent of $\hat{Q}^{n} \left( 0 \right)$ or converge to a deterministic limit), so we will focus on the convergence of $\hat{X}^{n}$.

Recall from Eq.\ \eqref{eq:Xhatdef} and Eq.\ \eqref{eq:Xhatinputdef} that $\hat{X}^{n}$ is built up from six different types of one-dimensional processes, namely \begin{equation}
\label{proc}n^{1 - \beta} \bar{X}^{n}_{1,k},\:\:n^{1 - \beta} \bar{X}^{n}_{2,k},\:\:n^{1 - \beta} \bar{X}^{n}_{3,k},\:\:n^{1 - \beta} \bar{X}^{n}_{4,kl},\:\:n^{1 - \beta} \bar{X}^{n}_{5,kl},\:\:\mbox{and}\:\: n^{1 - \beta} \bar{X}^{n}_{6,kl}.\end{equation}
Clearly, joint weak convergence of all these processes implies weak convergence of $\hat{X}^{n}$, as a result of the CMT.

To establish joint weak convergence of these processes, we will exploit the fact that some of them are independent or converge to deterministic functions. For those processes that are independent or converge to a deterministic function, joint convergence follows automatically (cf.\ \cite[Th.\ 11.4.4]{whitt2002} and \cite[Th.\ 11.4.5]{whitt2002}).
We will prove joint weak convergence of the six processes in (\ref{proc})  in four steps. We will then summarize these steps in a fifth and final step.

\emph{Step 1.}
Some of the six types of processes (namely the the processes $n^{1 - \beta} \bar{X}^{n}_{2,k}$ and $n^{1 - \beta} \bar{X}^{n}_{5,kl}$) converge in probability to deterministic functions. Indeed,
\begin{align*}
n^{1 - \beta} \bar{X}^{n}_{2,k} \left( t \right) &= n^{1 - \beta} \left( \frac{1}{n} \lambda^{n}_{k} - \lambda_{k} \right)^{\transpose} \int_{0}^{t} K^{n} \left( s \right) \, \dd s,
\end{align*}
which converges uoc in probability to the deterministic process 
\begin{align*}
\hat{\lambda}_{k}^{\transpose} \int_{0}^{t} \pi \, \dd s = \int_{0}^{t} \hat{\lambda}_{k}^{\pi} \, \dd s
\end{align*}
by Eq.\ \eqref{eq:lambdahatdef} and Theorem \ref{thm:occupationmeasureconvergence}. The process $n^{1 - \beta} \bar{X}^{n}_{5,kl}$ satisfies
\begin{align*}
n^{1 - \beta} \bar{X}^{n}_{5,kl} \left( t \right)
&= n^{1 - \beta} \left( \frac{1}{n} \mu^{n}_{kl} - \mu_{kl} \right)^{\transpose} \int_{0}^{t} \frac{1}{n} Q^{n}_{k} \left( s \right) K^{n} \left( s \right) \, \dd s \\
&= n^{1 - \beta} \left( \frac{1}{n} \mu^{n}_{kl} - \mu_{kl} \right)^{\transpose} \int_{0}^{t} \left( \frac{1}{n} Q^{n}_{k} \left( s \right) - \rho_{k} \left( s \right) \right) K^{n} \left( s \right) \, \dd s \\
&\phantom{{} = {}} {} + n^{1 - \beta} \left( \frac{1}{n} \mu^{n}_{kl} - \mu_{kl} \right)^{\transpose} \int_{0}^{t} \rho_{k} \left( s \right) K^{n} \left( s \right) \, \dd s.
\end{align*}
Consider the last two terms above. The first of these converges uoc in probability to the zero process. This follows immediately from the fluid limit in Lemma~\ref{lem:fluidlimit} and the fact that $K^{n}$ is bounded by~$1$. The convergence of the second term, on the other hand, follows from similar arguments as the convergence of $n^{1 - \beta} \bar{X}^{n}_{2,k}$ by combining Eq.\ \eqref{eq:muhatdef} and Theorem \ref{thm:occupationmeasureconvergence}. In this case, the limiting process is given by
\begin{align*}
\int_{0}^{t} \mu^{\pi}_{kl} \rho_{k} \left( s \right) \, \dd s.
\end{align*}
Consequently, $n^{1 - \beta} \bar{X}^{n}_{5,kl}$ converges uoc in probability to this process, too.

\emph{Step 2.}
We have reduced our task to proving joint convergence of $n^{1 - \beta} \bar{X}^{n}_{1,k}$, $n^{1 - \beta} \bar{X}^{n}_{3,k}$, $n^{1 - \beta} \bar{X}^{n}_{4,kl}$, and $n^{1 - \beta} \bar{X}^{n}_{6,kl}$.
A closer inspection of these processes reveals that each of them comes in one of two flavors. On the one hand, $n^{1 - \beta} \bar{X}^{n}_{3,k}$ and $n^{1 - \beta} \bar{X}^{n}_{6,kl}$ can be written as a stochastic integral with respect to the semimartingale $G^{n}$. On the other hand, $n^{1 - \beta} \bar{X}^{n}_{1,k}$ and $n^{1 - \beta} \bar{X}^{n}_{4,kl}$ can be viewed as a random time change of the Poisson martingales
\begin{align*}
n^{1 - \beta} \bar{A}^{n}_{k} \left( t \right) = n^{1 - \beta} \left( \frac{1}{n} A_{k} \left( n t \right) - t \right)
,\:\:\:\:\:\:
n^{1 - \beta} \bar{S}^{n}_{kl} \left( t \right) = n^{1 - \beta} \left( \frac{1}{n} S_{kl} \left( n t \right) - t \right),
\end{align*}
respectively. Indeed, defining the random times
\begin{align*}
\tau^{n}_{1,k} \left( t \right)
= \int_{0}^{t} \frac{1}{n} \lambda_{k}^{n} \left( J^{n} \left( s \right) \right) \, \dd s,\:\:\:\:\:\:\tau^{n}_{4,kl} \left( t \right)
= \int_{0}^{t} \mu_{kl}^{n} \left( J^{n} \left( s \right) \right) \frac{1}{n} Q_{k}^{n} \left( s \right) \, \dd s,
\end{align*}
we see that $n^{1 - \beta} \bar{X}^{n}_{1,k} = n^{1 - \beta} \bar{A}^{n}_{k} \circ \tau^{n}_{1,k}$ and $n^{1 - \beta} \bar{X}^{n}_{4,kl} = n^{1 - \beta} \bar{S}^{n}_{kl} \circ \tau^{n}_{4,kl}$.

The idea is to prove joint weak convergence of $n^{1 - \beta} \bar{X}^{n}_{3,k}$ and $n^{1 - \beta} \bar{X}^{n}_{6,kl}$, together with $n^{1 - \beta} \bar{A}^{n}_{k}$, $n^{1 - \beta} \bar{S}^{n}_{kl}$, and the random times $\tau^{n}_{1,k}$ and $\tau^{n}_{4,kl}$. Say that they have limits $X_{3,k}$, $X_{6,kl}$, $B_{k}$, $B_{kl}$, $\tau_{1,k}$, and $\tau_{4,kl}$, respectively. Then, under mild conditions, \cite[Th.~13.2.2]{whitt2002} guarantees the joint weak convergence of $n^{1 - \beta} \bar{X}^{n}_{3,k}$, $n^{1 - \beta} \bar{X}^{n}_{6,kl}$, $n^{1 - \beta} \bar{A}^{n}_{k} \circ \tau^{n}_{1,k}$, and $n^{1 - \beta} \bar{S}^{n}_{kl} \circ \tau^{n}_{4,kl}$.
Additionally, \cite[Th.~13.2.2]{whitt2002} identifies the corresponding limits as $X_{3,k}$, $X_{6,kl}$, $B_{k} \circ \tau_{1,k}$, and $B_{kl} \circ \tau_{4,kl}$. The conditions of this theorem will be clearly met in our case, as all limits will turn out to be continuous.

\emph{Step 3.}
Thus, before applying \cite[Th.~13.2.2]{whitt2002}, we have to prove joint weak convergence of $n^{1 - \beta} \bar{X}^{n}_{3,k}$ and $n^{1 - \beta} \bar{X}^{n}_{6,kl}$, together with $n^{1 - \beta} \bar{A}^{n}_{k}$, $n^{1 - \beta} \bar{S}^{n}_{kl}$, and the random times $\tau^{n}_{1,k}$ and $\tau^{n}_{4,kl}$. To prove this joint weak convergence, observe that the processes $n^{1 - \beta} \bar{X}^{n}_{3,k}$ and $n^{1 - \beta} \bar{X}^{n}_{6,kl}$ are independent from the Poisson martingales, because the background process is independent from these martingales by assumption. Consequently, $n^{1 - \beta} \bar{X}^{n}_{3,k}$, $n^{1 - \beta} \bar{X}^{n}_{6,kl}$, and the Poisson martingales converge jointly as soon as $n^{1 - \beta} \bar{X}^{n}_{3,k}$ and $n^{1 - \beta} \bar{X}^{n}_{6,kl}$ converge jointly and the Poisson martingales converge jointly.

However, in general $n^{1 - \beta} \bar{X}^{n}_{3,k}$, $n^{1 - \beta} \bar{X}^{n}_{6,kl}$, and the Poisson martingales are not independent from $\tau^{n}_{1,k}$ and $\tau^{n}_{4,kl}$. Nevertheless, dealing with the random time changes is relatively straightforward: they converge uoc in probability to deterministic functions, which is a sufficient condition for joint weak convergence. Note that, similar to $n^{1 - \beta} \bar{X}^{n}_{2,k}$ and $n^{1 - \beta} \bar{X}^{n}_{5,kl}$, the process $\tau^{n}_{1,k}$ converges to the continuous function $\tau_{1,k} \left( t \right) = \int_{0}^{t} \lambda^{\pi}_{k} \, \dd s$ and the process $\tau^{n}_{4,kl}$ converges to the continuous function $\tau_{4,kl} \left( t \right) = \int_{0}^{t} \mu^{\pi}_{kl} \rho_{k} \left( s \right) \, \dd s$. 

Now consider the Poisson martingales and recall that they are all independent. For ease of exposition, take $n^{1-\beta} \bar{A}^{n}_{k} = n^{1/2 - \beta} \sqrt{n} \bar{A}^{n}_{k}$. We know from the MCLT that $\sqrt{n} \bar{A}^{n}_{k}$ converges weakly to a standard Brownian motion $B_{k}$ (cf.\ \cite[Th.\ 4.2]{PTW2007}). Since $1/2 - \beta = \min \left\lbrace 0 , \left( \alpha - 1 \right) / 2 \right\rbrace$, it follows that $n^{1-\beta} \bar{A}^{n}_{k}$ converges weakly to $\ind{\alpha \geq 1} B_{k}$. Similarly, $n^{1 - \beta} \bar{S}^{n}_{kl}$ converges weakly to $\ind{\alpha \geq 1} B_{kl}$, where $B_{kl}$ is a standard Brownian motion. Because of the Poisson martingales being independent, joint convergence holds as well and the Brownian motions $B_{k}$ and $B_{kl}$ are independent.

\emph{Step 4.}
It remains to show joint weak convergence of the processes $n^{1 - \beta} \bar{X}^{n}_{3,k}$ and $n^{1 - \beta} \bar{X}^{n}_{6,kl}$. We observed earlier that these processes are stochastic integrals with respect to the continuous semimartingale $G^{n}$. More concretely, we have
\begin{align*}
n^{1 - \beta} \bar{X}^{n}_{3,k} \left( t \right)
&= \int_{0}^{t} n^{1 - \beta - \alpha/2} \lambda^{\transpose}_{k} \, \dd n^{\alpha/2} G^{n} \left( s \right)
\end{align*}
and
\begin{align*}
n^{1 - \beta} \bar{X}^{n}_{6,kl} \left( t \right)
&= \int_{0}^{t} n^{1 - \beta - \alpha/2} \rho_{k} \left( s- \right) \mu^{\transpose}_{kl} \, \dd n^{\alpha/2} G^{n} \left( s \right).
\end{align*}
The form of these integrals suggests that we may apply Theorem \ref{thm:markovintegralconvergence}.

Let us check that the conditions of Theorem \ref{thm:markovintegralconvergence} are satisfied here. First, recall the definition of $\beta$ in Eq.\ \eqref{eq:betadef}. In particular, it holds that $1 - \beta - \alpha/2 = 0$ for $\alpha \leq 1$, $1 - \beta - \alpha/2 < 0$ for $\alpha > 1$, and $1 - \beta - \alpha < 0$.
Next, observe that $n^{1 - \beta - \alpha/2} \lambda^{\transpose}_{k}$ converges to $\ind{\alpha \leq 1} \lambda^{\transpose}_{k}$ uoc. In addition, $n^{1 - \beta - \alpha/2} \rho_{k} \mu^{\transpose}_{kl}$ converges to $\ind{\alpha \leq 1} \rho_{k} \mu^{\transpose}_{kl}$ uoc. Third, notice that $\rho_{k}$ is a function of bounded variation and that $\rho_{k}$ does not depend on $n$. Consequently, both $n^{- \alpha/2} n^{1 - \beta - \alpha/2} \lambda^{\transpose}_{k}$ and $n^{- \alpha/2} n^{1 - \beta - \alpha/2} \rho_{k} \mu^{\transpose}_{kl}$ are finite-variation processes whose total variation processes converge to the zero process uoc.

Combining all these observations, we see that the conditions of Theorem \ref{thm:markovintegralconvergence} are satisfied. Then Theorem \ref{thm:markovintegralconvergence} implies the joint weak convergence of the stochastic integrals $n^{1 - \beta} \bar{X}^{n}_{3,k}$ and $n^{1 - \beta} \bar{X}^{n}_{6,kl}$ to \[\int_{0}^{t} \ind{\alpha \leq 1} \lambda^{\transpose}_{k} \, \dd B_{\mathrm{b}} \left( s \right)\:\:\mbox{ and}\:\:\ \int_{0}^{t} \ind{\alpha \leq 1} \rho_{k} \mu^{\transpose}_{kl} \, \dd B_{\mathrm{b}} \left( s \right),\] respectively.

\emph{Step 5.}
Summarizing, we have shown joint weak convergence of the eight processes 
\begin{align*}
n^{1 - \beta} \bar{X}^{n}_{2,k},\:\: n^{1 - \beta} \bar{X}^{n}_{5,kl}, \:\:n^{1 - \beta} \bar{A}^{n}_{k},\:\: n^{1 - \beta} \bar{S}^{n}_{kl}, \:\:\tau^{n}_{1,k},\:\: \tau^{n}_{4,kl}, \:\:n^{1 - \beta} \bar{X}^{n}_{3,k}, \:\:\text{ and }\:\: n^{1 - \beta} \bar{X}^{n}_{6,kl}
\end{align*}
to, respectively,
\begin{align*}
&\int_{0}^{t} \hat{\lambda}_{k}^{\pi} \, \dd s, \:\:\int_{0}^{t} \mu^{\pi}_{kl} \rho_{k} \left( s \right) \, \dd s, \:\:\ind{\alpha \geq 1} B_{k},\:\: \ind{\alpha \geq 1} B_{kl}, \:\:\tau_{1,k}, \:\:\tau_{4,kl}, \\
&\int_{0}^{t} \ind{\alpha \leq 1} \lambda^{\transpose}_{k} \, \dd B_{\mathrm{b}} \left( s \right), \text{ and } \int_{0}^{t} \ind{\alpha \leq 1} \rho_{k} \mu^{\transpose}_{kl} \, \dd B_{\mathrm{b}} \left( s \right).
\end{align*}
This implies joint weak convergence of the six processes
\begin{align*}
n^{1 - \beta} \bar{X}^{n}_{2,k}, \:\:n^{1 - \beta} \bar{X}^{n}_{5,kl}, \:\:n^{1 - \beta} \bar{X}^{n}_{1,k}, \:\:n^{1 - \beta} \bar{X}^{n}_{4,kl}, \:\:n^{1 - \beta} \bar{X}^{n}_{3,k},\:\: \text{ and }\:\: n^{1 - \beta} \bar{X}^{n}_{6,kl}
\end{align*}
to, respectively,
\begin{align*}
&\int_{0}^{t} \hat{\lambda}_{k}^{\pi} \, \dd s,\:\: \int_{0}^{t} \mu^{\pi}_{kl} \rho_{k} \left( s \right) \, \dd s, \:\:\ind{\alpha \geq 1} B_{k} \circ \tau_{1,k},\:\: \ind{\alpha \geq 1} B_{kl} \circ \tau_{4,kl}, \\
&\int_{0}^{t} \ind{\alpha \leq 1} \lambda^{\transpose}_{k} \, \dd B_{\mathrm{b}} \left( s \right), \text{ and } \int_{0}^{t} \ind{\alpha \leq 1} \rho_{k} \mu^{\transpose}_{kl} \, \dd B_{\mathrm{b}} \left( s \right).
\end{align*}
Now applying the CMT, and using that $\ind{\alpha \geq 1} B_{k} \circ \tau_{1,k}$ and $\ind{\alpha \geq 1} B_{kl} \circ \tau_{4,kl}$ have the same distributional properties as the integrals \[\ind{\alpha \geq 1} \int_{0}^{t} \sqrt{ \lambda_{k}^{\pi} } \, \dd B_{1,k} \left( s \right)\:\:\:\mbox{and}\:\:\ind{\alpha \geq 1} \int_{0}^{t} \sqrt{ \mu_{kl}^{\pi} \rho_{k} \left( s \right) } \, \dd B_{kl} \left( s \right),\]
the statement of Theorem \ref{thm:diffusionlimit} follows immediately.

\bibliographystyle{plain}
\bibliography{ref_networks}

\end{document}